\documentclass{amsart}    

\usepackage{graphicx,import}

\usepackage{amssymb}
\usepackage{amsmath} 
\usepackage{amscd} 
\usepackage[all]{xypic}
\usepackage{url}
\usepackage{color}

\newcommand{\D}{\mathcal{D}} 
\newcommand{\va}{\varphi}

\newcommand{\fk}{\mathfrak k}

\newcommand{\fn}{\mathfrak n}

\newcommand{\Fl}{{\rm Fl}}

\newcommand{\R}{\mathbb{R}}

\newtheorem{prop*}{Proposition}

\newtheorem{thm*}{Theorem}

\newtheorem{lemma*}{Lemma}

\newtheorem{cor*}{Corollary}

\newtheorem{rem*}{Remark}

\newtheorem{def*}{Definition}

\newtheorem{exam}{Example}

\begin{document}

\title[Local control of a  mechanism with the growth vector $(4,7)$]{Local geometric control of a certain mechanism with the growth vector $(4,7)$
}


\thanks{}

\maketitle

\begin{center}
    \author{\bf Jaroslav Hrdina\,$^a$, Lenka Zalabov\'a\,$^b$}
	
	\vspace{6pt}
	
	\small
	
	$^a$\,
	Institute of Mathematics, \\
	Faculty of Mechanical Engineering,  Brno University of Technology,\\
	Technick\' a 2896/2, 616 69 Brno, Czech Republic,
\\	hrdina@fme.vutbr.cz \\
	
	$^b$\,
	Institute of Mathematics, 
	Faculty of Science, University of South Bohemia, \\
	Brani\v sovsk\' a 1760, \v Cesk\' e Bud\v ejovice, 370 05, Czech Republic, 
	and \\
	Department of Mathematics and Statistics, 
	 Faculty of Science, Masaryk University, \\
	  Kotl\' a\v rsk\' a 2, Brno, 611 37, Czech Republic,	\\ lzalabova@gmail.com
	
\end{center} 

\thispagestyle{empty}

\vspace{7pt}
%%%%%%%%%%%%%%%%%%%%%%%%%%%%%%%%%%%%%%%%%%%%%%%%%%%%%%%%%%%%
%%%%%%%%%%%%%%%%%%%%%%%%%%%%%%%%%%%%%%%%%%%%%%%%%%%%%%%%%%%%
\begin{abstract}
	We study local control of the mechanism with the growth vector $(4,7)$. We study controllability and extremal trajectories on the nilpotent approximation as an example of the control theory on Lie group. We give solutions of the system an show examples of local extremal trajectories. 
	\end{abstract}

\section{Introduction}

Originally,  the general trident snake robot has been introduced in \cite{I04}. 
Let us recall that the trident robot is a mechanism composed of three snake--legs, each connected to an equilateral triangle root block in its vertices \cite{hvnm2,I04,I10,PT} for further details.
Generally, the branches can be multi--link, assumed that each link has its own passive wheel, which provides footing for the robot. Active elements, which affect controllability, are placed on branches. 
Its  simplest non--trivial version, corresponding to one--links,  has  been  mainly discussed, see  e.g. \cite{I10,hvnm2,PT}.
In this case, the control distribution is that of the growth vector $(3,6)$ \cite{mya1}.

We are interested in the modification corresponding to one or more prismatic joints such that the control distribution will be a that of the growth vector $(4,7)$.
Local controllability of such robot is given by the appropriate Pfaff system of ODEs.  
The solution space gives a control system $\dot q = \sum u_i X_i$ where the vector fields $X_1, X_2, X_3, X_4$ describe the horizontal distribution and $u: \mathbb R \to \mathbb R^4$ is the control of the system.
Consequently, the system is controllable by Chow--Rashevsky theorem \cite{ABB,J,cg09}, see Section \ref{s2}. 

We construct a nilpotent approximation to get nilpotent Lie algebra and corresponding Lie group to study the controllability of approximated left invariant control system, see Section \ref{s3}. We study geometric properties and symmetries of the nilpotent approximation in Section \ref{s4}.
We use the theory of Hamiltonians and Pontryagin's  maximum principle to study local control and extremal trajectories, see Section \ref{s5}. In particular, we provide analysis of the system and present explicit solutions. 
%We demonstrate a periodic input to show some local behavior of our mechanisms and we discuss convenient parameters of the input in Section \ref{s6}. 

\section{Analysis of the mechanisms} \label{s2}
In this Section we describe a 
mechanism that is a modification of the trident snake robot
(for more details see \cite{hvnm2,I04,I10,PT}).
We present our mechanisms as a new example of non-holonomic
systems with multi--generators and discuss local controllability of our mechanism  based on the principles of non-holonomic mechanics.

\subsection{Description of the mechanism and its movement}
\begin{figure}[ht]
	\centering 
	\includegraphics[scale=0.25]{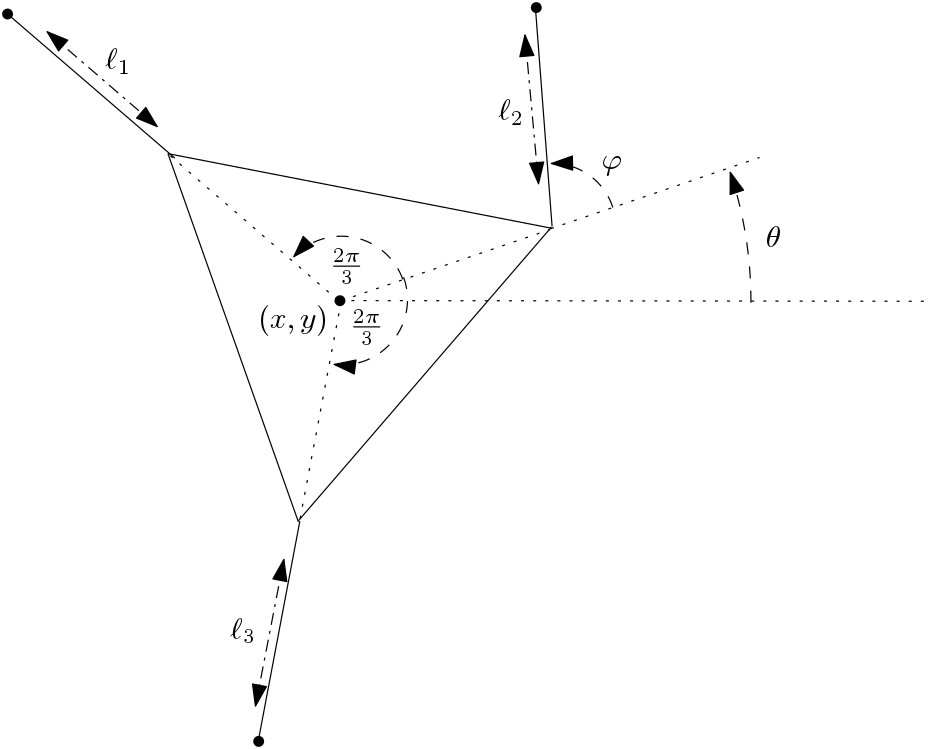}
	\caption{Description of the mechanism}
	\label{t}
\end{figure} 

In the sequel, we study a generalized trident mechanism which consists of a root block in the shape of an equilateral triangle with three $1$--link branches
with variable length.  Each of the branch is connected to one vertex of the root block and they form prismatic joints.
Second leg is in addition connected
	to the root block via the revolute joint, connection of the two remaining is fixed. Each link has a
	passive wheel on its branches, which is assumed neither to slip,
	nor slide sideways.

The configuration space of the planar mechanism in question corresponds to a manifold that locally coincides with $\mathbb R^7$ (but not globally). Since we study local problems, we can consider the configuration space to be $\R^7$ with the coordinates  $(x, y, \theta, \varphi, \ell_1, \ell_2, \ell_3)$. With respect to a fixed coordinate system, the first three coordinates describe completely the global position of the mechanism in the plane. The coordinates $x$ and $y$ give the position of the centre of mass of the root block in $\R^2$ and $\theta$ gives the amount of the counter-clockwise rotation.
% from the positive $x$--axis to the positive $y$--axis. 
Remaining four coordinates represent the input for the mechanism's active elements. Thus as active elements we consider the revolute joint of the branch $\ell_2$ with the root block, namely $\varphi$, and prismatic joints, which can change absolute lengths of branches $\ell_1, \ell_2$ and $\ell_3$, see Figure \ref{t}.

To provide the description of the robot's movement, we shall only point out that we do not get any singularities as long as the robot's configuration satisfies that $0 < \ell_i$ and $\varphi $ is not too far from $0$.

Using the method of moving frame, the kinematic model is the
set of equations  
for $i$--th wheel position in the form
\begin{align} 
\begin{split}
x_i &= x + \cos (\theta + \alpha_i) + \ell_i \cos (\theta + \alpha_i),  \label{ro1} \\
y_i &= y + \sin (\theta + \alpha_i) + \ell_i \sin (\theta + \alpha_i ), 
\end{split}
\end{align} 
for $i=1,3$, where $\alpha_1 = - \frac{2}{3} \pi, \: \alpha_3 = \frac{2}{3} \pi$, and
\begin{align}
\begin{split}
x_2 &= x + \cos (\theta ) + \ell_2 \cos (\theta + \va), \label{ro2} \\
y_2 &= y + \sin (\theta) + \ell_2 \sin (\theta  + \va). 
\end{split}
\end{align}

Consequently, we differentiate the position equations \eqref{ro1} and  \eqref{ro2} with respect to time $t$ and obtain the  velocity equations as follows
\begin{align*}
\dot x_i &= \dot x - \sin (\theta + \alpha_i)\dot \theta - \ell_i \sin (\theta + \alpha_i )\dot \theta  
+ \cos (\theta + \alpha_i ) \dot \ell_i, \\
\dot y_i &= \dot y + \cos (\theta + \alpha_i)\dot \theta + \ell_i \cos (\theta + \alpha_i)\dot \theta
+ \sin (\theta + \alpha_i) \dot \ell_i,
\end{align*}
for $i=1,3$, where $\alpha_1 = - \frac{2}{3} \pi, \: \alpha_3 = \frac{2}{3} \pi$, and
\begin{align*}
\dot x_2 &= \dot x - \sin (\theta )\dot \theta - \ell_2 \sin (\theta + \va)(\dot \theta  + \dot \va)
+ \cos (\theta  + \va) \dot \ell_2, \\
\dot y_2 &= \dot y + \cos (\theta )\dot \theta + \ell_2 \cos (\theta  + \va)(\dot \theta  + \dot \va)
+ \sin (\theta  + \va) \dot \ell_2.
\end{align*}
The conditions preventing slipping lead to 
the non-holonomic constraints of the form
\begin{align*}
0 &= (-\sin (\theta + \alpha_i), \cos (\theta + \alpha_i))
\cdot (\dot x_i, \dot y_i ),  \\
0 &= (-\sin (\theta + \va), \cos (\theta  + \va))
\cdot (\dot x_2, \dot y_2 ) ,
\end{align*}
where $i=1,3$, $\alpha_1 = - \frac{2}{3} \pi, \: \alpha_3 = \frac{2}{3} \pi$  
and $\cdot$ is the Riemannian scalar product on the Euclidean space $\mathbb R^2$. 
We obtain the following differential kinematic equations,
so the movement of the mechanism is described by the Pfaff system of three nonlinear homogeneous equations 
\begin{align} 
\label{pf} 
\begin{split}
0 & = - \sin (\theta -\frac{2 \pi}{3}   ) \text{dx} + 
\cos (\theta -\frac{2 \pi}{3} ) \text{dy} + (1+ \ell_1) \text{d}\theta  ,\\
0 & =- \sin (\theta  + \va) \text{dx} + \cos (\theta  + \va) \text{dy} + (-\cos(\varphi)  + \ell_2 ) \text{d}\theta -\ell_2 \text{d} \va  , \\
0 & = - \sin (\theta + \frac{2 \pi}{3}) \text{dx} + \cos (\theta + \frac{2 \pi}{3} ) \text{dy}+  (1 + \ell_3) \text{d} \theta.
\end{split}
\end{align}

\subsection{Local controllability of the system}
\label{loccon}
The space of solutions to the system \eqref{pf} forms four dimensional distribution on the configuration space, so--called horizontal distribution. 
It follows from \eqref{pf} that the solution space always contains the vector fields  
$X_2:=\partial_{\ell_1}$, $X_3:=\partial_{\ell_2}$ and $X_4:=\partial_{\ell_3}$ as generators.
In the case $\ell_2 \neq 0$ the last generating vector field $X_1$ is a combination of $\partial_x, \partial_y, \partial_\theta$ and $\partial_\varphi$ which is generically complicated and we do not need to write it here in the biggest generality while $X_1= \partial_{\va}$ in the case $\ell_2=0$. 
In fact, from the mechanical point of view, zero length of each leg makes no sense, so we suppose in the following that $\ell_i>0$.

We can equivalently rewrite the solution space of the Pfaff  system \eqref{pf} in the following form

\begin{align}  \label{pff}
\dot q = u_1 X_1 (q) + u_2 X_2 (q) + u_3 X_3 (q) +u_4 X_4 (q),
\end{align} 
where $  q = ( x,  y,  \theta,  \va,  \ell_1,  \ell_2, \ell_3)$. This is a $4$--input
symmetric  affine  
control system. 
In general, controllability  of symmetric affine systems is completely characterized by the controllability Lie algebra by
Chow--Rashevsky's theorem \cite{J, S,MZS}. 
Control system \eqref{pff} satisfies the Chow's  condition 
at the point $q$ if 
$Lie(X_1, X_2, X_3, X_4)(q) = T_q \mathbb R^7,$
where the controllability Lie algebra $Lie(X_1, X_2, X_3, X_4)$ is the Lie algebra generated by $X_1,X_2,X_3,X_4$. In this case point $q$ is regular 
and can be connected with any point in a suitable
neighbourhood of $q$ by a horizontal trajectory.
If the Chow's condition is satisfied at all points (of a connected space) then any two points can be joined by a horizontal trajectory 
and the system is locally controllable.

In our case, mechanical description leads to an observation that local controllability depends on the shape of the mechanism only, not on its configuration in the plane. 
In other words, regular points
have to be invariant with respect to the rigid body (Euclidean)  transformations of the plane. 
So we choose 
$x=y=0$ and $\theta={\pi \over 2}$ 
without loss of generality.
In particular, for points of the form $q_0=(0,0,{\pi \over 2},\varphi,\ell_1,\ell_2,\ell_3)$
the vector fields 
\begin{align*}
X_1 &= \partial_x+\frac{(\ell_1-\ell_3)\sqrt{3}}{3L} \partial_{y}
-\frac{1}{L} \partial_{\theta}
+\frac{\sin(\varphi)\sqrt{3}(  \ell_1-  \ell_3)+3\cos(\varphi)(  L+1 )+3 \ell_2}{3 \ell_2 L} \partial_{\varphi} \\
X_2 &= \partial_{\ell_1}, \: \:\: 
X_3 = \partial_{\ell_2} ,\:\:\:
X_4 =  \partial_{\ell_3}, 
\end{align*}
generate the solution space of our Pfaff system, where we denote $L=\ell_1+\ell_3+2$.
Moreover, at $q_0$, the controllability Lie algebra is obtained by the Lie bracket operation as
\begin{align*}
X_{12}&:= [X_1,X_2] = 
 \frac{-2(\ell_3+1 )}{\sqrt{3} L^2} \partial_{y}
- \frac{1}{L^2} \partial_{\theta} 
 + \frac{ -2   \sin(\varphi)( \ell_3 +1 )+\sqrt{3} \cos(\varphi)+\sqrt{3} \ell_2}{ \sqrt{3} \ell_2 L^2} \partial_{\varphi},
  \\
X_{13}&:=[X_1,X_3] =\frac{ \sin(\varphi)( \ell_1- \ell_3)+ \sqrt{3}  \cos(\varphi)( L+1) }{\sqrt{3}  \ell_2^2 L} \partial_{\varphi} , \\
X_{14}&:=[X_1,X_4] =  
\frac{2\ell_1+2}{\sqrt{3} L^2} \partial_{y}
- \frac{1}{L^2} \partial_{\theta}
+\frac{ 2  \sin(\varphi) (\ell_1+1)+\sqrt{3} \cos(\varphi)+\sqrt{3} \ell_2}{\sqrt{3} \ell_2 L^2} \partial_{\varphi}, \\
\end{align*}
and remaining brackets are trivial.
Then the  matrix $\bar G$ consisting of coordinates of vector fields $X_1,X_2,X_3,X_4,X_{12},X_{13},X_{14}$
spans full tangent space $\mathbb R^7$
as long as $\text{det} ( \bar{G}(q_0) )\neq 0$ and the 
system is locally controllable at $q_0.$
Thus our system is locally controllable at $q_0$ (and without loss of generality everywhere because of mechanical meaning) because  $\ell_2 \neq 0$ and $L \neq 0$.

We consider the filtration  $\Delta^{1} = \langle X_1,\dots, X_4 \rangle \subset \Delta^{2} =
\langle X_1,\dots, X_4, X_{12},X_{13},X_{14} \rangle
$. We have  $\dim \Delta^1 (q) =4$ and $\dim \Delta^2 (q) =7$ at all points 
and we have filtration with the growth vector $(4,7)$.

\subsection{Remark on  corresponding dynamical systems}
If we restrict our considerations only to nontrivial movements 
of the root block, i.e. movements in the  $X_1$ direction or in the direction
of iterated bracket of $X_1$ and $X_i$ for $i=2, 3, 4$, 
we can use methods of dynamic pairs.

In general, each control affine system
\begin{align*}
 \dot x = X(x) + \sum_{j=1}^m  u_j Y_j (x)
\end{align*}
 on a manifold $M$, where $X, Y_1, \dots, Y_m$ are smooth vector fields on $M$ and 
 $u =(u_1, \dots , u_m)^T$
are controls, defines a dynamic pair $(X, V)$, where $V$ is a distribution spanned by $Y_1, . . . , Y_m$.
Then there is a sequence of distributions defined inductively, using Lie bracket, by
$ \label{filtration}
V^0:= V, \: \: \: V^{i+1}:=V^i+[X,V^i]
$
and one imposes the regularity conditions \cite{jk13,doze}:
\begin{align*} 
 \text{rk} \: V^i = (i+1)m, \; \: \text{ for }i=0,\dots,k, \\
  V^k \oplus \langle X \rangle = TM \label{R2}.
\end{align*} 

Then our control system \eqref{pff} can be adapted to this situation as $X:=f(q) X_1$, where $f$ is an arbitrary non--zero function and $Y_1:=X_2,Y_2:=X_3,Y_3:=X_4$, and
the regularity conditions are satisfied. Indeed,   
\begin{align*}
V^0 &= \langle Y_1,Y_2,Y_3 \rangle ,\\ 
V^1 &= \langle Y_1,Y_2,Y_3 ,[X,Y_1],[X,Y_2],[X,Y_3]\rangle , 
\end{align*}
such that  $\text{rk} \: V^0 = 3=m$, $\text{rk} \: V^1 = 6=2m$ 
 and $V^1 \oplus \langle X \rangle = TM $. In fact, each $f(q)$ defines a dynamical system with specific drift and all of them are regular.

 \subsection{Remark on dual curvature} 	Following \cite{M,DZ}, curvature of a distribution $H$ on a manifold $Q$ is a linear bundle map  
 	$ F: \wedge^2 H \to TQ / H  $
 	defined by 
 	$ F(X,Y) = - [X,Y] \text{ mod } H. $
 	Denote by $H^{\bot}$  the bundle of covectors that annihilate $H$. Since the curvature $F$ is a linear bundle map, the dual of curvature is a linear map 
 	$ F^* :H^{\bot} \to  \wedge^2 H^* ,$
 	 called the dual curvature. 
 	Because our distribution is equipped with the growth vector $(4,7)$, the space $\text{Im} (F^*)$
 	is three--dimensional subspace
 	of  $\wedge^2 H^* $.
 	We can define the Pfaffian 
 	$H^{\bot}  \to \wedge^4 H^*  $
 	as $ \mu \mapsto F^* (\mu ) \wedge F^* (\mu )$. One can see that 
 	$\wedge^4 H_q^*$, $q \in Q,$ 
 	is one--dimensional vector space and the Pfaffian may be understood as a real valued quadratic form on $H^{\bot}$  by choosing a volume form. 
 	Then 
 	possible signatures for the Pfaffian are
 	$ (3,0), (2,1),(2,0), (1,1), (1,0), (0,0). $ Note that the signatures $(p,r)$ and $(r,p)$ must be considered as identical because $\wedge^4 H_q^*$ is not oriented.
 
 From parabolic geometry point of view, the generic distribution with growth vector $(4,7)$ corresponds to quaternionic contact structures and split quaternionic contact structures for signature $(3,0)$ and $(2,1),$ respectively. The distribution corresponding to our mechanism has signature $(0,0)$ and we will see that it is a parabolic geometry called generalized path geometry \cite{CS}. Let us point out that there exist different modifications of trident snake robot that lead to the growth vector $(4,7)$. However, it turned out that all of them have non-regular signature. We have not found a mechanism with signature $(3,0)$ or $(2,1)$, yet.

\section{Nilpotent approximation}\label{s3}

We recall a constructive method to approximate vector fields of a nonlinear control system by a similar system on the same configuration space. The method leads to an approximate distribution which has a nilpotent basis. The techniques of nilpotent approximation  have been developed by 
various researchers, see e.g. \cite{as87,h86}.
We recall the following concept of orders of functions or vector fields and
distribution weights. Let $X_i$, $i=1, \dots, m$ denote a family of smooth vector fields on
a manifold $M$ and $C^{\infty}(p)$ the set of germs of smooth functions at $p \in M$.
For $f \in  C^{\infty}(p)$ we say that Lie derivatives $X_i f, X_i X_j f, \dots$ are non–holonomic derivatives of $f$ of order $1,2,\dots$. The non–holonomic derivative of order $0$ of $f$ at $p$ is $f(p)$.
Then the non–holonomic order $\text{ord}_p (f )$ of $f$ at $p$
is the biggest integer $k$ such that all non–holonomic derivatives of $f$ of order smaller than $k$ vanish at $p$, i.e.
\begin{align*}
\text{ord}_p(f) = 
\text{min} \bigg 
\{
s \in \mathbb N: \exists i_1 \dots i_s \in
\{ 1,\dots ,m \} \: \text{s.t.}
( X_{i_1} \cdots X_{i_s} f) (p) \neq 0
\bigg \}.
\end{align*}
Denote by $VF(p)$ the set of germs of smooth vector fields at $p \in M$.
Then the notion of non–holonomic order extends to vector fields as follows:
For $X \in VF(p)$ the non–holonomic order $\text{ord}_p (X)$ of $X$ at $p$ 
is a real number defined by
\begin{align*}
\text{ord}_p(X) = 
\text{sup} \bigg 
\{
\sigma \in \mathbb R: 
\text{ord}_p(Xf) \geq 
\sigma + \text{ord}_p(f), 
\forall f \in 
C^{\infty}(p)
\bigg \}. 
\end{align*}
Note that $\text{ord}_p (X) \in  \mathbb Z$. Moreover, the zero vector field $X \equiv  0$ has infinite order, i.e. $\text{ord}_p (0) = \infty$. Furthermore, $X_1,\dots ,X_m$ are of order $\geq -1$, $[X_i , X_j ]$ of
order $\geq  -2$, etc.

\subsection{Construction}
We construct a 
nilpotent approximation of the 
distribution with respect to the given filtration at point 
$q_0= (0,0,\frac{\pi}{2},0,1,1,1) $. 
We use Bella\"{i}che algorithm, which may be found in \cite{b96}.
Let us point out that all constructions are local in the neighbourhood of $q_0$.
In our case, as the first step of Bella\"{i}che algorithm, 
the adapted frame 
\begin{align}
X_1, X_2, X_3, X_4, X_{12}, X_{13}, X_{14}
\end{align}
was chosen. 
Then we use four local coordinates $(x, \ell_1,\ell_2,\ell_3)$ as first four adapted coordinates.  The others can be obtained from the original coordinate system by an affine change in the form
\begin{align}
\begin{split} 
y_1 & = -2 x-2 \sqrt{3} y-8 \theta ,\\
y_2 & =\frac{4}{5} \va - \frac{4}{5} x+\frac{8}{5} \theta ,\\
y_3 & = -2 x+2 \sqrt{3} y-8 \theta,
\end{split}\label{trans}
\end{align}
where we use the following conditions:
\begin{align*}
\partial_x |_{q_0} &= X_1(q_0), \: \:
\partial_{\ell_1} |_{q_0} = X_2(q_0), \: \:
\partial_{\ell_2} |_{q_0} = X_3(q_0), \:\:
 \partial_{\ell_3} |_{q_0} = X_4(q_0), \\
 \partial_{y_1} |_{q_0} &= X_{12}(q_0), \: \:
 \partial_{y_2} |_{q_0} = X_{13}(q_0), \: \:
  \partial_{y_3} |_{q_0} = X_{23}(q_0).
\end{align*} 
Finally, note that the adapted coordinates $(x, \ell_1, \ell_2, \ell_3, y_1, y_2 , y_3 )$ are privileged coordinates in the case of adapted frame with growth vector $(4,7)$  \cite{J}.

Following \cite{J}, we group together the monomial vector fields in Taylor expansions of the same weighted degree and thus we express $X_i$ as a series
$$X_i=X_i^{(-1)}+X_i^{(0)}+X_i^{(1)}+\cdots,$$
for $i=1,2,3,4$, where $X_i^{(s)}$ is a homogeneous vector field of order $s$. 	
By \cite[Proposition 2.3.]{J} we set $N_i := X_i^{(-1)}$ for $i=1,2,3,4$. The family of vector fields $(N_1,N_2,N_3,N_4)$ forms called first order approximation of $(X_1,X_2,X_3, X_4)$ at $q_0$ and generates a nilpotent Lie algebra of step $r=2$, i.e. all brackets of length greater than 2 are zero.	
The family $(N_1, N_2, N_3 , N_4)$ is called the (homogeneous) nilpotent
approximation of $(X_1, X_2, X_3 , X_4)$ at $q_0$ associated with coordinates $(x, \ell_1, \ell_2, \ell_3, y_1, y_2 , y_3 )$.

In the sequel, by the above algorithm, we obtain the following vector fields:
\begin{align*}
N_1&= \partial_x -\left(\frac{-\sqrt{3}}{2}x+\ell_1-1\right) \partial_{y_1} -(\ell_2-1) 
\partial_{y_2} - \left(\frac{\sqrt{3}}{2}x +\ell_3-1\right)\partial_{y_3},
\\
N_2&= \partial_{\ell_1} ,\: \: \:
N_3 = \partial_{\ell_2}, \: \: \:
N_4 =\partial_{\ell_3}.
\end{align*}
In particular, the family of vector fields  $(N_1,N_2,N_3,N_4)$ is the nilpotent approximation of vector fields $(X_1, X_2,X_3,X_4)$ at $(0,0,\frac{\pi}{2},0,1,1,1)$ in the coordinates $(x,y,\theta,\varphi, \ell_1, \ell_2,\ell_3)$, while it is the point  
$(0,1,1,1,-4 \pi,\frac{4}{5} \pi,-4 \pi)$ in the coordinates
$(x,\ell_1,\ell_2,\ell_3, y_1, y_2,y_3)$.

The remaining three vector fields $N_{12},N_{13},N_{14}$ are generated by Lie brackets of $(N_1,N_2,N_3, N_4)$ as 
\begin{align*}
N_{12}&= [N_1,N_2] =  \partial_{y_1},\:\:\: 
N_{13}= [N_1,N_3] =\partial_{y_2},\:\:\:
N_{14}=[N_1,N_4] = \partial_{y_3}.
\end{align*}
Note that due to the linearity of coefficients of $(N_1,N_2,N_3, N_4)$, the coefficients of $(N_{12},N_{13},N_{14})$ must be constant.

To show how the nilpotent approximation affects integral curves of the 
distribution and the resulting control, we compute  the Lie brackets of relevant vector fields.
\begin{figure}[ht]
	\centering 
	\includegraphics[scale=0.09]{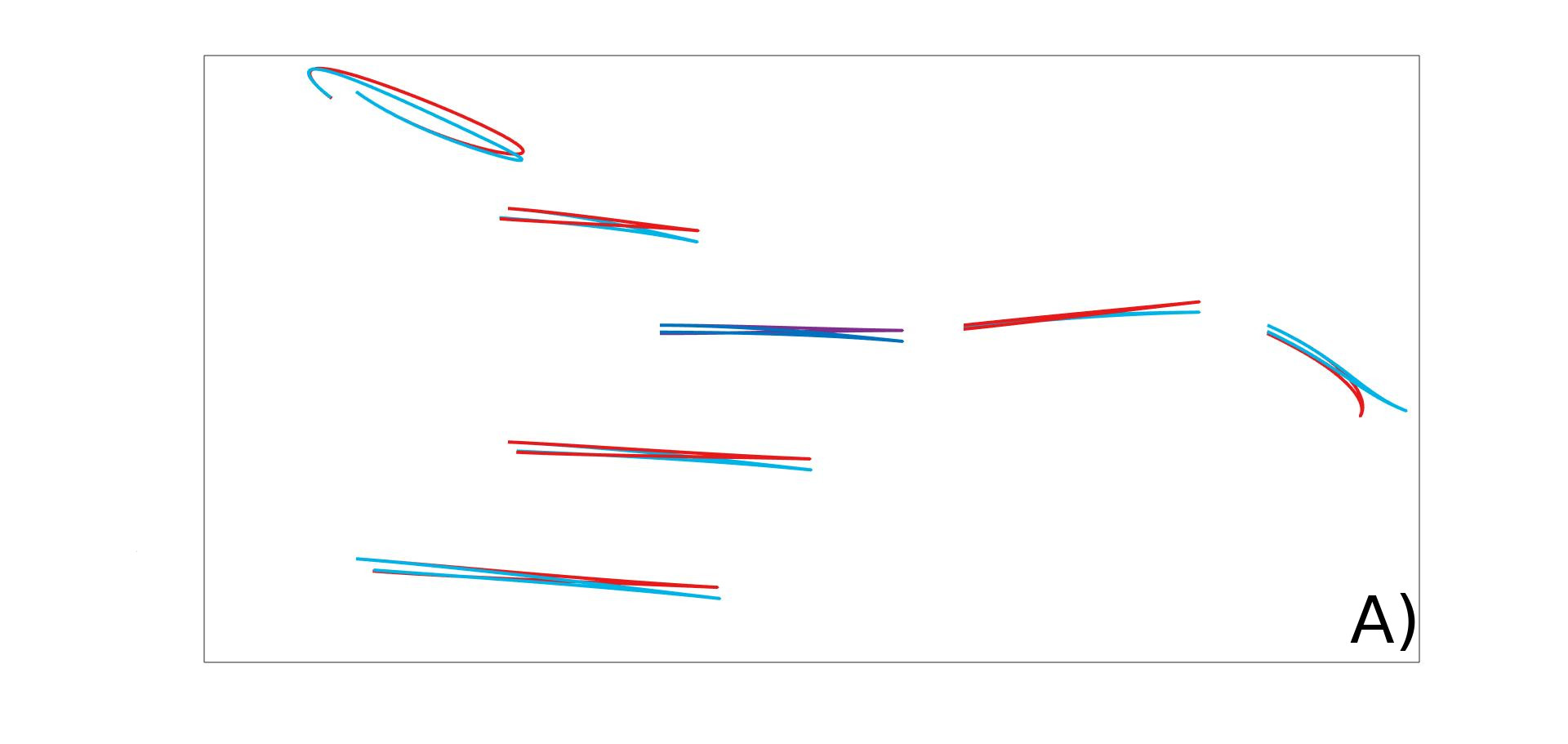}
	\includegraphics[scale=0.09]{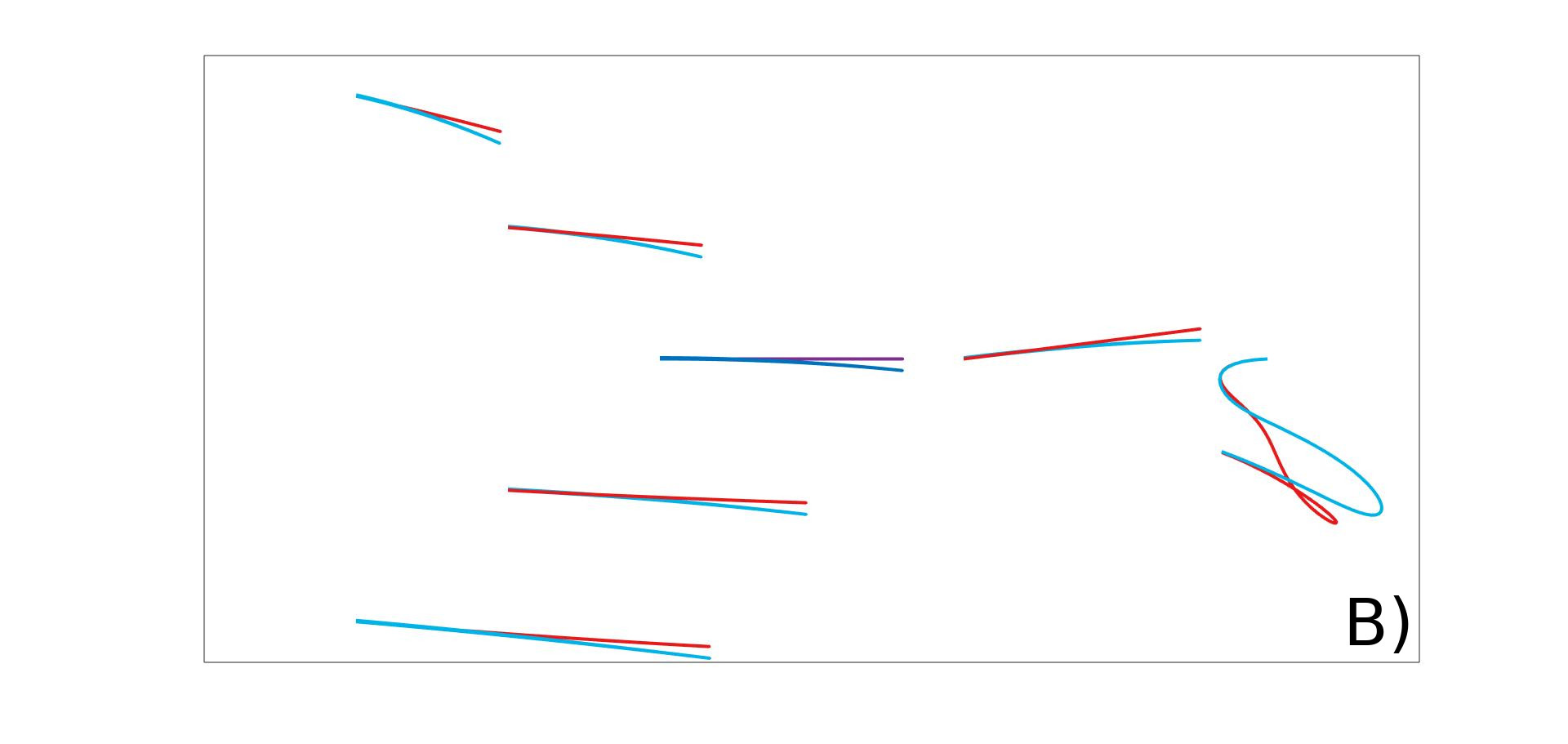}
	\includegraphics[scale=0.09]{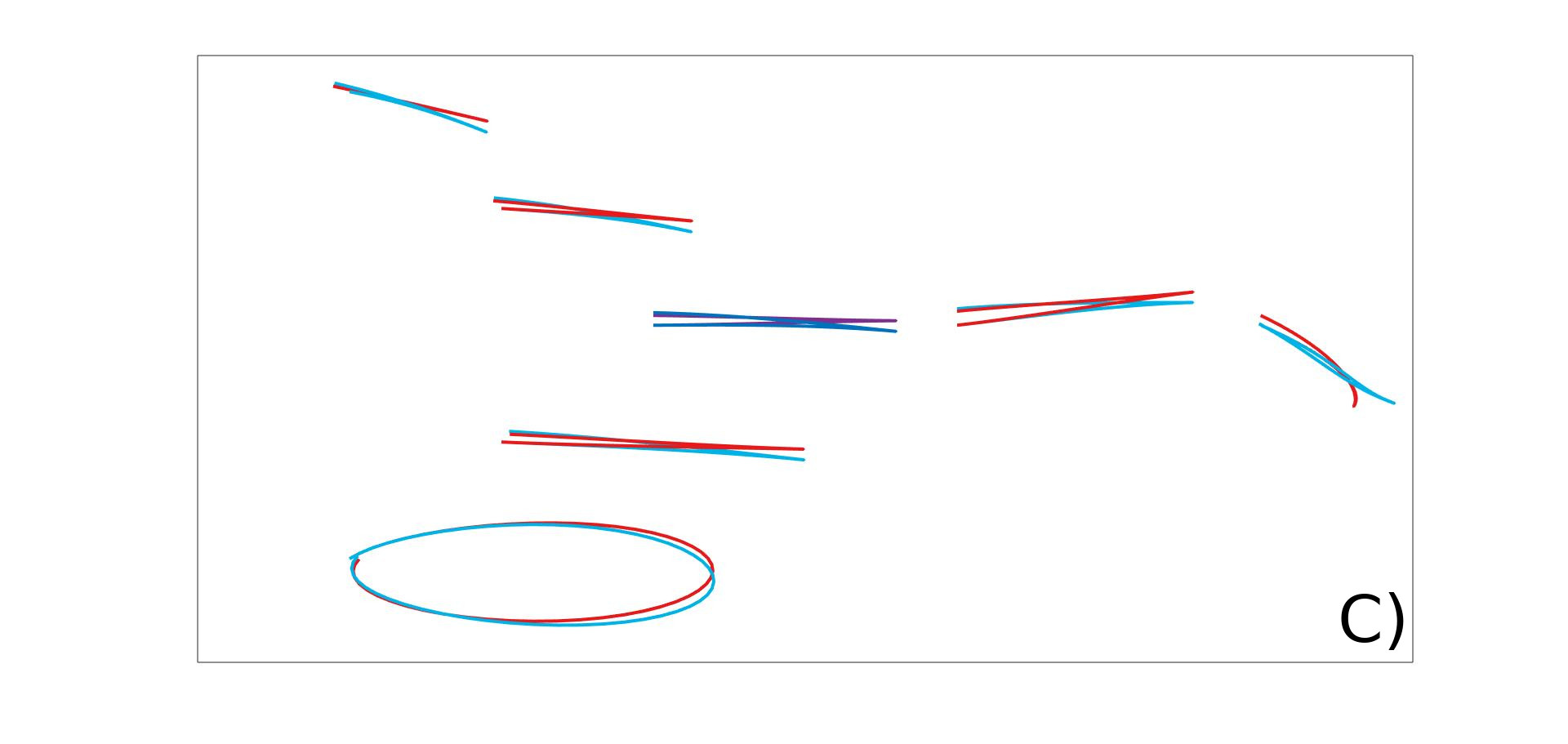}
    \includegraphics[scale=0.09]{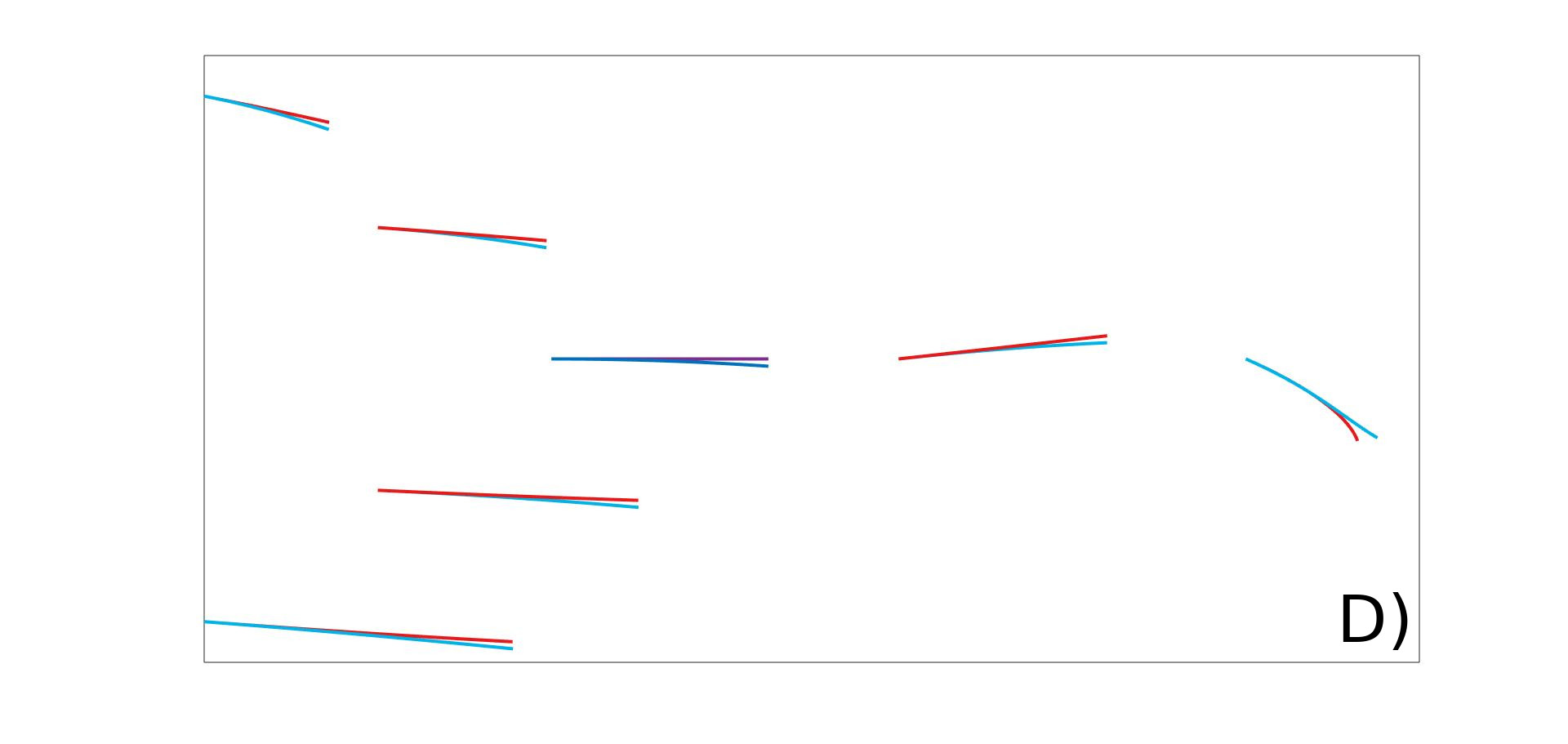}	
	\caption{ Motions of the mechanism in directions: A) $N_{12}$,  $X_{12}$, B) $N_{13}$,  $X_{13}$, C) $N_{14}$,  $X_{14}$, D) $N_{1}$,  $X_{1}$, }
	\label{nilp}
\end{figure}
In Fig. \ref{nilp}, there is a comparison of the Lie bracket motions in the original distribution (red line) and in the nilpotent approximation (blue line). The following figures show the trajectories of the root centre point, vertices and wheels when a particular Lie bracket motion is realized.
To simulate the bracket motion of the nilpotent approximation, we choose the initial state $q_0$ and apply the periodic input on couples $(N_1, N_2)$, $(N_1,N_3)$ and $(N_1,N_4)$ to receive the displacement approximately parallel to $[N_1,N_2]$, $[N_1,N_3]$ and $[N_1,N_4]$ respectively \cite{MZS}. More precisely, for one cycle we apply
	\begin{align*}
	u_1(t) & = -A \omega  \sin( \omega t), \:\:\:
	u_i(t)  = A \omega  \cos (\omega t ),\:\:\:
	u_j(t)  = 0 ,\:\:\:
	u_k(t)  = 0 
\end{align*}
	for $i \in \{2,3,4 \}$ and $j,k \in \{ 2,3,4\} - \{ i\}$, respectively, all with respect to the control system \eqref{control-system}, amplitude $A=0.4$ and angular speed $ \omega= \frac{2 \pi }{50}$. Then we apply the same process to original vector fields $X_1,X_2,X_3$ and $X_4$.

\subsection{Control theory on nilpotent Lie groups}
\label{section3.2}

The above construction led to vector 
fields $N_1,N_2,N_3,N_4,N_{12},N_{13},N_{14}$ which 
determine a $7$--dimensional nilpotent Lie algebra $\fn$. 
The computation gives that 
\begin{align*}
[N_1,N_2] &= N_{12}, \:\:\:
[N_1,N_3] = N_{13}, \:\:\:
[N_1,N_4] = N_{14}
\end{align*}
and the remaining brackets are trivial.
The corresponding connected simply connected nilpotent Lie group $N\simeq \R^7$ with the nilpotent Lie algebra $\fn$ is then endowed, in new coordinates $(x,\ell_1,\ell_2,\ell_3,y_1,y_2,y_3),$ with the following group structure
\begin{align}
\begin{pmatrix}
 x \\  \ell_1 \\ \ell_2 \\ \ell_3 \\  y_1 \\ y_2 \\ y_3 
\end{pmatrix} \times 
\begin{pmatrix}  
 \bar x \\  \bar \ell_1 \\ \bar \ell_2 \\ \bar \ell_3 \\ \bar y_1 \\ \bar y_2 \\  \bar y_3 \\
\end{pmatrix}
=
\begin{pmatrix}  
x + \bar x \\ \ell_1 + \bar \ell_1 \\ \ell_2 + \bar \ell_2 \\ \ell_3 + \bar \ell_3 \\ y_1+\bar y_1 +\frac{\sqrt{3}}{2} x \bar{x}   - \ell_1 \bar{x}  \\ y_2+ \bar y_2 - \ell_2\bar{x}  \\ y_3+ \bar y_3  - \frac{\sqrt{3}}{2} x \bar x - \ell_3\bar{x} \\
\end{pmatrix} \label{ngroup}
\end{align}
and the vector fields $N_1,$ $N_2,$ $N_3,$ $N_4,$ $N_{12},$ $N_{13},$ $N_{14}$ are left--invariant with respect to the left action given by the group structure.
In particular, the vector fields $N_i$ for $i=1,2,3,4$ determine a left--invariant distribution $\D$ on $N$, which has the growth vector $(4,7)$ everywhere. 

Altogether, $\D \subset TN$ defines an invariant $4$--input symmetric  affine 
control system
\begin{align} \label{control-nil}
\dot q = u_1N_1+u_2N_2+u_3N_3+u_4N_4,
\end{align} 
with $q=(x,\ell_1,\ell_2,\ell_3,y_1,y_2,y_3)$, that approximates the original control system. It clearly satisfies the Chow's condition and is controllable.

\section{Infinitesimal symmetries} \label{prvni-cast}\label{s4}
We focus on basic geometric properties and infinitesimal symmetries of the nilpotent approximation. By infinitesimal symmetries we mean vector fields such that their flows preserve the geometric structure at any time \cite{s09}. In our case, infinitesimal symmetries preserve the distribution and also the metric.
\subsection{Generalized path geometries and their symmetries}
\label{parabolic}
Let us discuss the geometric structure that occurs behind the control problem \eqref{control-nil}. 
In the previous section, we constructed a nilpotent Lie group $N$ with a filtered nilpotent Lie algebra $\fn$ with the growth vector $(4,7)$, where the $4$--dimensional distribution $\D$ is generated by the left--invariant fields $N_1$, $N_2$, $N_3$ and $N_4$. Consider subbundles $E=\langle N_1\rangle$, $V=\langle N_2, N_3, N_4\rangle$ in $TN$. One can see from the structure of Lie brackets that the following holds:
\begin{enumerate}
	\item $E \cap V =0$,
	\item Lie bracket of two sections of $V$ is a section of $E \oplus V$, and
	\item for sections $\xi \in \Gamma(E)$ and $\nu \in \Gamma(V)$ and a point $q \in N$, the equation $[\xi,\nu] (q) \in E_q \oplus V_q$ implies that $\xi(q)=0$ or $\nu(q)=0$.
\end{enumerate}    
Such geometric structures are usually called  \emph{(generalized) path geometries} (in dimension $7$) \cite[Section 4.4.3]{CS}. 

General theory \cite{CS} says that generalized path geometries have finite--dimen\-sional Lie algebras of symmetries and in the case of $7$--dimensional manifolds, the maximal possible dimension is $24$. This is the case of generalized path geometries that are locally equivalent to a generalized flag manifold $PSL(5,\R)/P_{1,2}$, where by $PSL(5,\R)$ we denote the projectivised special linear group with Lie algebra $\frak{sl}(5,\R)$, and by $P_{1,2}$ the stabiliser of the flag of a line in a plane  for the projectivised standard action of $PSL(5,\R)$. In particular, the Lie algebra of symmetries for such generalized path geometry is exactly the simple Lie algebra $\frak{sl}(5,\R)$ and the symmetries with a fixed point form its $17$--dimensional parabolic subalgebra.

There is a general method to find all infinitesimal symmetries of a nilpotent filtered structure \cite{s09}. One can apply this method to our structure $(N,E \oplus V \subset TN)$ and it turns out that the symmetry algebra is of dimension $24$ and is exactly $\frak{sl}(5,\R)$. In other words, symmetries of $(N,E \oplus V \in TN)$ are left multiplications by elements of a Lie group with the Lie algebra $\frak{sl}(5,\R)$ and the structure is left--invariant with respect to this action. From this point of view, $N \simeq PSL(P,\R)/P_{1,2}$.
Altogether, nilpotent approximation forms a flat generalized path geometry.

\begin{rem*}
The concept of Cartan geometries \cite{Sharpe} generalizes the concept of Klein geometries \cite{kl1,kl2} and generalized flag manifolds are special cases of Klein geometries for the case of parabolic subgroups in semisimple groups. The Cartan's generalization then leads to a wide theory of parabolic geometries \cite{CS}, that are curved versions of flag manifolds. Generalized path geometries are examples of such parabolic geometries \cite[Definition  4.4.3.]{CS}. 
\end{rem*}

\subsection{Sub--Riemannian structure and its symmetries}
To study extremal trajectories in the next section, we need the sub--Riemannian structure on the nilpotent approximation. We consider a control metric $g$ in $\D=\langle N_1, N_2,N_3,N_4 \rangle$ such that the fields $N_i$ for $i=1,2,3,4$ are orthogonal and have the length one with respect to $g$.  This clearly determines a left--invariant sub--Riemannian structure $g$ of $\D$ (with respect to the action given by group structure \eqref{ngroup} on $N$).

Let us now focus on the symmetries of the nilpotent control problem $(M,\D=E+V,g)$. Thus we are interested in symmetries that preserve not only the flat generalized path geometry, but also the control metric. The symmetry algebra $\fk$ of $(M,E+V,g)$ clearly is a subalgebra of $\frak{sl}(5, \R)$. In fact, both $\fk$ and  $\frak{sl}(5, \R)$ contain the same nilpotent subalgebra that reflects the nilpotent group structure and that acts effectively and transitively on $N$. 
It is generated by vector fields
\begin{align*}
w_1 &:=-  \partial_x-
\frac{\sqrt{3}}{2}  \partial_{\ell_1} +\frac{\sqrt{3}}{2} x \partial_{y_3},\\
w_2 &:= \partial_{\ell_1}-x \partial_{y_1}, \\
w_3 &:= \partial_{\ell_2}-x \partial_{y_2}, \\
w_4 &:= \partial_{\ell_3}-x \partial_{y_3} ,\\
w_{12} &:= \partial_{y_1}, \\
w_{13} &:= \partial_{y_2}, \\
w_{14} &:= \partial_{y_3}.
\end{align*}

In general, $\fk$ also contains symmetries preserving arbitrary fixed point. Since all points are equivalent, we can fix the origin $o=(0,0,0,0,0,0,0)$.
 If we study symmetries preserving the origin, intuition suggests that the sub--Riemannian metric $g$ shall be preserved by an orthogonal algebra $\frak{so}(4)$. However, each such symmetry shall also preserve the control distribution $\D$ and its decomposition into $E$ and $V$. Thus it acts trivially on the $1$--dimensional subspace $E_o$ and restricts to the action of $\frak{so}(3)$ on the $3$--dimensional subspace $V_o$. Direct computation gives that there really is the symmetry algebra $\frak{so}(3) \subset \fk$ preserving the origin generated by fields
\begin{align}
\begin{split}
&v_1:=-\ell_3\partial_{\ell_2}+\ell_2\partial_{\ell_3}-({\sqrt{3}x^2 \over 4} -x+y_3)\partial_{y_2}-(x-y_2)\partial_{y_3},\\
&v_2:=\ell_3\partial_{\ell_1}-\ell_1\partial_{\ell_3}+({\sqrt{3}x^2 \over 4} -x+y_3)\partial_{y_1}
+({\sqrt{3}x^2 \over 4}+x-y_1)\partial_{y_3},\\
&v_3:= -\ell_2 \partial_{\ell_1}+\ell_1 \partial_{\ell_2}+ (x-y_2) \partial_{y_1} -({\sqrt{3}x^2 \over 4}+x-y_1)\partial_{y_2}.
\end{split}
\end{align}	
We can write it also in a `vector--matrix--like' notation as 
\begin{align*}
\begin{pmatrix}
v_1 \\ v_2 \\ v_3
\end{pmatrix}&=
\begin{pmatrix}
0 & -\ell_3 & \ell_2 \\
 \ell_3 & 0 & -\ell_1\\
 -\ell_2 & \ell_1 & 0
\end{pmatrix}
\begin{pmatrix}
\partial_{\ell_1} \\ \partial_{\ell_2} \\ \partial_{\ell_3}
\end{pmatrix}\\ &+
\begin{pmatrix}
0 & -({\sqrt{3}x^2 \over 4} -x+y_3) & -(x-y_2) \\
{\sqrt{3}x^2 \over 4} -x+y_3 & 0 & {\sqrt{3}x^2 \over 4} +x-y_1\\
 x-y_2 & -({\sqrt{3}x^2 \over 4} +x-y_1) & 0
\end{pmatrix}
\begin{pmatrix}
\partial_{y_1} \\ \partial_{y_2} \\ \partial_{y_3}
\end{pmatrix}.
\end{align*}
One can verify directly that $\mathcal{L}_{v_i}V \subset V$ and $\mathcal{L}_{v_i}N_1= 0$ and $\mathcal{L}_{v_i}g =0$ for $i=1,2,3$, where $\mathcal{L}$ denotes the Lie derivative, and that $[v_1,v_2]=-v_3$, $[v_1,v_3]=v_2$ and $[v_2,v_3]=-v_1$.

\subsection{Properties of the $\frak{so}(3)$--action}
The action of $\frak{so}(3)$ 
is clear if the action of the generators $v_1, v_2, v_3$ is explained: the action decomposes into two independent actions first of which is on $\R^3$ given by $\partial_{\ell_1}$, $\partial_{\ell_2}$ and $\partial_{\ell_3}$ and the second one on $\R^4$ given by $\partial_x$, $\partial_{y_1}$, $\partial_{y_2}$ and $\partial_{y_3}$. Moreover, the algebra $\frak{so}(3)$ acts trivially on $\partial_x$, which then defines a one--dimensional invariant subspace. The action of the algebra $\frak{so}(3)$ on $(\partial_{\ell_1},\partial_{\ell_2},\partial_{\ell_3}) \in \R^3$ corresponds to the rotations around all axes passing through $(0,0,0)$ and this does not depend on the remaining variables $x$, $y_1$, $y_2$ and $y_3$. In particular, $v_i$ for $i=1,2,3$ corresponds to the rotation around the axis generated by $\partial_{\ell_i}$. Analogously, for arbitrary fixed $x$, the algebra $\frak{so}(3)$ acts on $(\partial_{y_1},\partial_{y_2},\partial_{y_3})$ via rotations around axes going through $(x+{\sqrt{3}x^2 \over 4}, 
x, x-{\sqrt{3}x^2 \over 4}).$ 
In particular, $v_i$ for $i=1,2,3$ corresponds to the rotation around the axis generated by $\partial_{y_i}$.

One can also see from the shape of the generators  that the rotation around $\partial_{\ell_i}$ is tied to the rotation around $\partial_{y_i}$ for $i=1,2,3$. 
Thus the action on $(\partial_{\ell_1},\partial_{\ell_2},\partial_{\ell_3})$ determines the action on $(\partial_{y_1},\partial_{y_2},\partial_{y_3})$ and vice versa.
In fact, $x$ parametrizes (in the coordinates $(x,\ell_1,\ell_2,\ell_3, y_1, y_2,y_3)$) the curve $(x,0,0,0,x+{\sqrt{3}x^2 \over 4} ,x , x -{\sqrt{3}x^2 \over 4}) \subset N$ which can be viewed as a curve of centres of the above `double--rotations'.

In particular, the nilpotent sub--Riemannian structure $(N,\D=V+E,g)$ is invariant with respect to the action, and we can study its action on curves passing through the origin. Assume that $c(t)$ is a (parametrized) curve such that $c(0)=o$. 
Consider the flow $\Fl^t_{v}$ of the infinitesimal symmetry $v:=a_1 v_1+a_2 v_2+a_3v_3 \in \frak{so}(3)$ for some $a_1, a_2, a_3$.  It clearly preserves the origin $o$. 
Assume that the point $c(t_0) \neq o$ for some $t_0 \neq 0$ is preserved by the action of $\Fl^t_{v}$. Then either the curve $c(t)$ is preserved by this action on $[0,t_0]$, or the action determines a family of curves of the same length from $o$ to $c(t_0)$ on $[0,t_0]$. 
In particular, if such a curve $c(t)$ that is not invariant with respect to the action of $\Fl^s_{v}$ is an extremal curve for the invariant control system, then it is no more minimiser after it reaches the point.
Moreover, if one finds one such point then its orbit with respect to the action of $\frak{so}(3)$ consists of such points. 
Indeed, a family of curves from $o$ to $c(t_0)$ is mapped to the  
family of curves of the same length from $o$ to $\hat c(t_0)$, where $\hat c(t_0)$ 
is image of $c(t)$ with respect to the action of $\mathfrak{so}(3)$.
%\end{color}

We can describe explicitly the set of such points that are fixed for the action of (the flow of) some infinitesimal symmetry $a_1v_1+a_2 v_2+a_3v_3$. First, one can check that each point of the curve $(x,0,0,0,x+{\sqrt{3}x^2 \over 4} ,x , x -{\sqrt{3}x^2 \over 4})$ is preserved by each such symmetry. Then the fixed points of any symmetry are given by axes of the corresponding `double--rotations'. Explicitly, the fixed points of the symmetry $a_1 v_1+a_2 v_2+a_3v_3$ form the set
$$
\left\{ (x,ka_1,ka_2,ka_3,x+{\sqrt{3}x^2 \over 4}+ka_1,x+ka_2 , x -{\sqrt{3}x^2 \over 4}+ka_3) : k \in \R \right\}.
$$

Let us finally say that for each $k_i$, $i=2,3,4$ it holds that $[N_1,k_2N_2+k_3N_3+k_4N_4]=k_2N_{12}+k_3N_{13}+k_4N_{14}$. So the triple 
$(N_1,k_2N_2+k_3N_3+k_4N_4,k_2N_{12}+k_3N_{13}+k_4N_{14})$ determines a subalgebra which has the structure of the Heisenberg algebra. One can see from the above mentioned that the action of the symmetry algebra $\frak{so}(3)$ simply maps each such Heisenberg subalgebra to another Heisenberg subalgebra.

\section{Pontryagin's maximum principle} \label{s5}
We study local control of particular mechanisms. We use Hamiltonian formalism and Pontryagin's maximum principle to find local length minimisers. We study the corresponding Hamiltonian system of ODEs in detail. Finally, we model several explicit movements of the mechanism. 

\subsection{Formulation of the problem}
Consider two configurations $q_1, q_2$ in the nilpotent approximation $N$.  Among all admissible curves $c(t)$, i.e. locally Lipschitz curves such that $c(0)=q_1$ and $c(T)=q_2$ that are tangent to $\D$ for almost all $t \in [0,T]$, we would like to find length minimisers with respect to $g$.

We would like to minimize the length $l$ among all the horizontal curves $c$, where the length is given by $l(c)=\int_0^T\sqrt{g(\dot c, \dot c)} dt$ for the control metric $g$. Let us recall that the distance between two points $q_1,q_2 \in N$ is defined as $d: M \times M \to [ 0, \infty]$,
$d (q_1,q_2) = \inf_{\{c \in \mathcal S_{q_1,q_2}\}} l(c),$
where 
$
\mathcal S_{q_1,q_2} = \{c : c(0)=q_1, c(T) =q_2, c \text{\ \ admissible} \}
$  \cite{ABB,cg09,J}.
However, since minimizing of the energy of a curve implies minimizing of its length, we will rather minimize energy of curves.

We know from Chow--Rashevsky theorem that the control system of the nilpotent approximation is controllable, see Section \ref{loccon}. In particular, any two points can be joined by a horizontal curve and the distance of arbitrary two points is finite \cite{ABB}. 

We study the following nilpotent control problem 
\begin{align} \label{control-system}
{d \over dt}q=u_1 
\left(
\begin{array}{c}
1\\
0\\
0\\
0\\
{\sqrt{3} \over 2}x-\ell_1+1\\
-\ell_2+1\\
-{\sqrt{3} \over 2}x-\ell_3+1\\
\end{array}
\right)
+u_2 
\left(
\begin{array}{c}
0\\
1\\
0\\
0\\
0\\
0\\
0\\
\end{array}
\right)
+u_3 
\left(
\begin{array}{c}
0\\
0\\
1\\
0\\
0\\
0\\
0\\
\end{array}
\right)
+u_4 
\left(
\begin{array}{c}
0\\
0\\
0\\
1\\
0\\
0\\
0
\end{array}
\right)
\end{align} 
for $q \in N$ and the control $u=(u_1,u_2,u_3,u_4) \in \R^4$ with the boundary condition $q(0)=(x_0,\ell_{10},\ell_{20},\ell_{30},y_{10},y_{20},y_{30})$ and $q(T)=(x_1,\ell_{11},\ell_{21},\ell_{31},y_{11},y_{21},y_{31})$ arbitrary but fixed, and we minimize
\begin{align} \label{min}
{1 \over 2}\int_{0}^{T} (u_1^2+u_2^2+u_3^2+u_4^2)\ dt.
\end{align}
Without loss of generality, we choose the origin $o=(0,0,0,0,0,0,0)$ as the initial point $q(0)=q_1$. Since we solve nilpotent control problem, we get curves starting at different points using the left action coming from the multiplication in $N$, see Section \ref{section3.2}.

\subsection{Hamiltonian formalism} 
Let us consider a cotangent bundle $T^*N \to N$ and the coordinate functions $h_i=\langle \lambda, N_i\rangle$, where $\lambda \in T^*N$. Then we consider the Hamiltonian of the maximum principle 
$$ 
H(\lambda,\nu)=u_1h_1+u_2h_2+u_3h_3+u_4h_4+{\nu \over 2} (u_1^2+u_2^2+u_3^2+u_4^2),$$
which is a family of smooth functions affine on fibres that is parametrized by controls $(u_1,u_2,u_3,u_4) \in \R^4$ and a real number $\nu \leq 0$. The Pontryagin's maximum principle can be formulated as follows \cite{ABB}:
Assume $(\bar u(t), c(t))$ to be a pair such that 
$c(t)$ is a length minimizer for $\eqref{control-system},\eqref{min}$ with control function $u=\bar u(t)$. Then there exist a Lipschitz curve $\lambda(t) \in T^*_{c(t)}M$ and a number $\nu \leq 0$ such that 
$$(\lambda(t),\nu) \neq 0, \:\: \dot \lambda(t)=\vec H_{\bar u(t)}(\lambda(t))$$ for Hamiltonian vector field $\vec H$ corresponding to $H$ and $H_{\bar u(t)}=\textrm{max} \ H(\lambda(t),\nu)$. 

If $\lambda(t)$ satisfies the principle for $\nu=0$, then it is called abnormal and it is called normal otherwise.
The abnormal minimiser is strictly abnormal, if it is not normal. It follows from the Goh condition \cite{sard, ABB} that there are no strictly abnormal minimisers in the case of $2$--step distributions.  
Thus we will focus on the case $\nu<0$ and we can normalize it in such a way that $\nu=-1$. 

The extreme is achieved when ${\partial H_u \over \partial u_i}=h_i-u_i=0$ for $i=1,2,3,4$ and this implies for the controls that $u_i=h_i$ for $i=1,2,3,4$. In this case, the Hamiltonian of the maximum principle is of the form
$$
H={1 \over 2}(h_1^2+h_2^2+h_3^2+h_4^2).
$$
Then the corresponding Hamiltonian system associated with $H$ 
is the following
\begin{align}
&\dot q
=h_1N_1(q)+h_2N_2(q)+h_3N_3(q)+h_4N_4(q),\label{hor} \\ 
&\dot h_i = \{H,h_i \}, \label{vert}
\end{align}
where $N_i$ are the generators of the sub--Riemannian structure, $q\in N$ and $\{ \ ,\ \}$ is the usual Lie--Poisson bracket \cite{ABB}.
Using the fact that $\{h_i,h_j \}=\langle \lambda,[N_i,N_j] \rangle$ we conclude that the fibre system $\eqref{vert}$ is of the form
\begin{align}
\begin{split}
\label{system-vert}
&{d \over dt}h_1(t)=-h_5(t)h_2(t)-h_6(t)h_3(t)-h_7(t)h_4(t)\\
&{d \over dt}h_2(t)=h_5(t)h_1(t) \\
&{d \over dt}h_3(t) = h_6(t)h_1(t)\\
&{d \over dt}h_4(t)= h_7(t)h_1(t)\\
&{d \over dt}h_5(t)={d \over dt}h_6(t)={d \over dt}h_7(t)=0.
\end{split}
\end{align}
The equations are clearly independent of the horizontal coordinates $q$. The base system has the form
\begin{align}
\begin{split}
\label{system-hor}
&{d \over dt}x(t)=h_1(t)\\
&{d \over dt}\ell_1(t)= h_2(t) \\
&{d \over dt}\ell_2(t) = h_3(t)\\
&{d \over dt}\ell_3(t) = h_4(t)\\
&{d \over dt}y_1(t)=\left(1+{\sqrt{3} \over 2}x(t)-\ell_1(t)\right) h_1(t)\\
&{d \over dt}y_2(t)=\left(1-\ell_2(t)\right) h_1(t)\\
&{d \over dt}y_3(t)=\left(1-{\sqrt{3} \over 2}x(t)-\ell_3(t)\right) h_1(t).\\
\end{split}
\end{align}

\subsection{Analysis of the fibre system}
Let us first discuss the fibre system \eqref{system-vert}. Obviously, functions $h_5$, $h_6$ and $h_7$ are constant. 
If they are all zero, then the functions $h_1$, $h_2$, $h_3$ and $h_4$ are constant, too. 

Let us denote the solution constants corresponding to $h_5$, $h_6$ and $h_7$ by $C_5$, $C_6$ and $C_7$, respectively, and assume that at least one of them is non--zero. Define $K \equiv \sqrt{C_5^2 + C_6^2 + C_7^2}$. Then we get
\begin{equation} \label{h1dd}
\ddot h_1 = -C_5 \dot h_2-C_6 \dot h_3-C_7 \dot h_4 = -(C_5^2 + C_6^2 + C_7^2)\,h_1 = -K^2 h_1.
\end{equation}
Since $K^2>0$, the solution of \eqref{h1dd} is
\begin{equation} \label{h1}
h_1 = C_{11} \cos(Kt) + C_{12} \sin(Kt),
\end{equation}
for some constants $C_{11}$ and $C_{12}$.
Now, if $C_5 \neq 0$, we have
\begin{equation*}
\dot h_2 = C_5 h_1 = C_5 (C_{11} \cos(Kt) + C_{12} \sin(Kt))
\end{equation*}
and hence
\begin{equation} \label{h2}
h_2 = \frac{C_5}{K} (C_{11} \sin(Kt) - C_{12} (\cos(Kt))+C_{13}.
\end{equation}
This analogously holds for $h_6$ and $h_7$ if $C_6 \neq 0$ and $C_7 \neq 0$, respectively. In the same way we get $h_3$ and $h_4$.
Indeed, only the equation for function $h_1$ merges everything together. 

We know from the above mentioned, that the functions $h_i$ for $i=1,2,3,4$ equal to the controls $u_i$ of the system. Thus we control each vector field $N_i$ with a function $u_i$ which is either a constant or which oscillates. In fact, $N_i$, $i=2,3,4$ corresponds to $\partial_{\ell_j}$, $j=1,2,3$ and reflects the movement of legs which is natural from the mechanical point of view. The field $N_1$ is the crucial field for the movement and if it is controlled by a constant, then the remaining fields $N_i$ are controlled by constant, too.

Let us point out that the choice of $h_5$, $h_6$ and $h_7$ to be zero or non--zero, in fact, corresponds to the choice of constants in the solution. Since that corresponds to the choice of initial conditions, we can interpret the choice of zero or non--zero solution $h_5$, $h_6$ and $h_7$  as the choice of initial conditions, equivalently.
Moreover, we are interested in the solution of the whole system $(\ref{system-vert},\ref{system-hor})$ and we are basically interested in the curves, that are images in $N$ with respect to the canonical projection  $T^*N \to N$. There can be curves of that type in $N$ that differ only by the parametrization. It is reasonable to consider the curves parametrized by the arc--length, only. These are exactly the images of restrictions of the canonical projection to the solutions of $(\ref{system-vert},\ref{system-hor})$  with initial conditions satisfying
$h_1^2+h_2^2+h_3^2+h_4^2=1$ at the initial point (and thus everywhere).

Finally, let us remark that the choice $h_6=h_7=0$ implies that $h_3, h_4$ are constants and the remaining equations result in the system $\dot h_1=-h_5h_2, \dot h_2=h_5h_1,\dot h_5=0$, which is the fibre system of the control problem on the Heisenberg group.

\subsection{Analysis of the base system}
Let us first say that it is enough to study solutions with the initial condition $q(0)=o$. Then one can use the action given by the multiplication on $N$ to get a solution starting at an arbitrary point in $N$.

The first four equations of the system depend  on $h_i$, $i=1,2,3,4$ only and can be computed by direct integration. Thus the equation $\dot x= h_1$ gives 
that either $x=C_1t$ for some constant $C_1$ in the case that $h_1$ is constant, or $$x={C_{11} \over K} \sin(Kt)-{C_{12} \over K}\cos(Kt)+{C_{12} \over K}$$ in the case that $h_1$ is of the form \eqref{h1}, where $K,C_{11},C_{12}$ are from the previous section. 
Then we compute directly from \eqref{h2} that 
\begin{align*}
\ell_1={C_5 \over K^2}(C_{11}-C_{11}\cos(Kt)-C_{12}\sin(Kt))+C_{13}t.
\end{align*}
Analogous observation can be made for $\ell_2$ and $\ell_3$.

From mechanical point of view, this simply means that each branch can either prolong or shorten or does not change its length in the first case, or oscillates in the second case. In particular, the choice $h_6=h_7=0$ reflects the situation when lengths $\ell_2$ and $\ell_3$ of the second and third branch are constant and the robot uses only the first branch of the length $\ell_1$. The same principles work for the remaining two branches and the corresponding choices of $h_i$. 

The equations for $y_i$ depend on $x$ and $\ell_i$ for $i=1,2,3,$ only. So we can find $y_i$ by considering closed subsystems for $x,\ell_i,y_i$ for $i=1,2,3$. Let us remark that $y_i$, $i=1,2,3$ have no evident mechanical meaning. One should use the transformation between these coordinates and the original ones to get some information about the behaviour of $y, \phi, \theta$. We do not mention all the solutions explicitly. We rather provide several examples to demonstrate the explicit paths and movements of the mechanism in the following section.

Let us finally remark that in the case that $h_5=h_6=h_7=0$, we can find for each solution an infinitesimal symmetry such that the solution is contained in the fixed points set of the symmetry. The situation is more complicated in the case when some of $h_5, h_6, h_7$ is non--zero.

\subsection{Examples of solutions}
Let us present several particular solutions (satisfying all initial conditions of the system and also constraints coming from the mechanical setting). To get an information about the movement of our mechanism, we transform the solutions into the original coordinates which express the kinematics of the mechanism. 
So we use the transformation between coordinates $(x,\ell_1,\ell_2,\ell_3,y_1,y_2,y_3)$ and $(x,y,\theta,\va,\ell_1,\ell_2,\ell_3)$ in the form
\begin{align*}
&\va = {5 \over 4}y_2+{3 \over 2}x+{1 \over 8}y_1+{1\over 8}y_3, \\
&\theta = -{1 \over 16}y_1-{1 \over 16}y_3-{1\over 4}x,\\
& y = -{1 \over 12}\sqrt{3}(y_1-y_3),
\end{align*}
which is the inverse transformation to the transformations \eqref{trans}.
We also add illustrative graphs of behaviour of the control parameters which are prismatic joints $\ell_i$ and the revolute joint $\varphi$.

\begin{exam}
In the case that $h_5=h_6=h_7=0$, we can choose $h_1 = {7\over 10},$ $h_2 = h_3 = {1 \over 2},$ $h_4 = {1 \over 10}$ and, with suitable choice of constants, we get the solution in the form 
\begin{align*}
&x = {7t \over 10},\ \ \  
y = ({7\sqrt{3} \over 600}-{49\over 800})t^2, \ \ \ \theta = {21t^2 \over 1600}-{21t \over 80}, \\ 
&\va = -{49 t^2 \over 200}+{21t \over 10}, \ \ \ 
\ell_1 =   
\ell_2 = {t \over 2},\ \ \  
\ell_3 = { t \over 10}.
\end{align*}
\begin{figure}[ht]
	\includegraphics[scale=0.28]{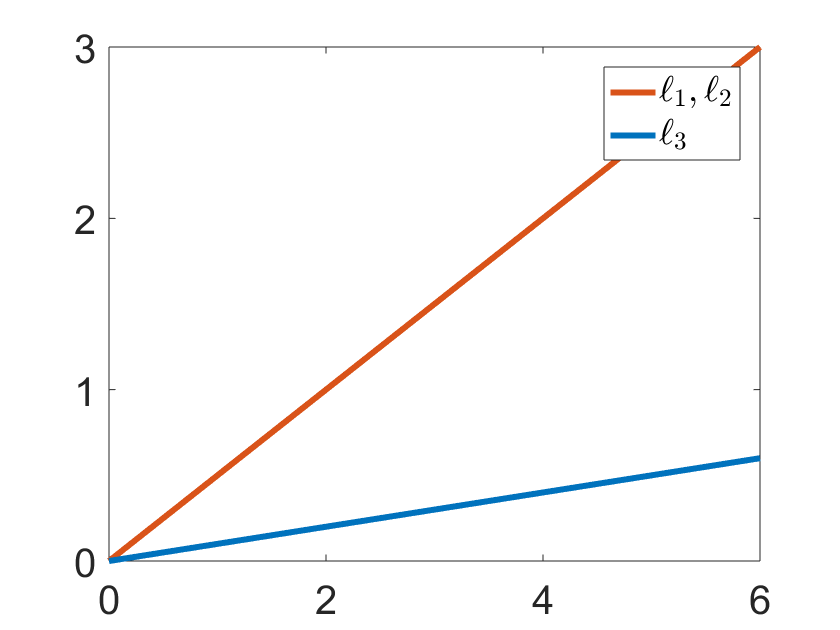}
	\includegraphics[scale=0.28]{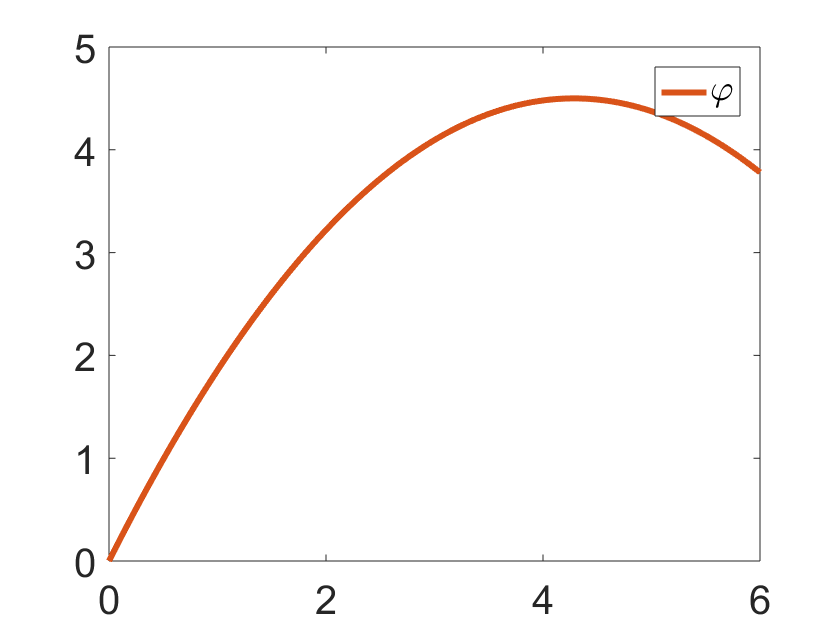}
	\caption{The graphs of control parameters: Example 1	}
\label{t1}
\end{figure} 
\end{exam}

\begin{exam}
In the case that $h_5 \neq 0$ and $h_6=h_7=0$, let us choose the constants in such a way that $h_5=1$ and $h_1 = -{1\over 2}\sin(t)+{1\over 2}\cos(t)$, $h_2 = {1 \over 2}\sin(t)+{1 \over 2}\cos(t)$ and $h_3 = h_4 =  {1 \over 2}$. Then we get the solution in the form

\begin{align*}
&x(t)={1\over 2}\sin(t)+{1\over 2}\cos(t)-{1\over 2},
\\&y(t)= -{\sqrt{3}\over 48}\left(\sqrt{3}(\sin(t)-1)(\cos(t)-1)+\cos(t)^2+t\cos(t)+(t-2)\sin(t)+t-1\right), 
\\&\theta(t)= -{1\over 64}\cos(t)^2+{t-10\over 64}\cos(t)+{t-12\over 64}\sin(t)-{t\over 64}+{11\over 64},
\\&\varphi(t)={1\over 32}\cos(t)^2+{36-11t\over 32}\cos(t)+{58 -11t\over 32}\sin(t)+{t\over 32}-{37\over 32}.
\\&  \ell_1(t)={1\over 2}\sin(t)-{1\over 2}\cos(t)+{1\over 2}, \ \ \ 
 \ell_2 = \ell_3 = {t\over 2}
\end{align*}

\begin{figure}[ht]
		\includegraphics[scale=0.28]{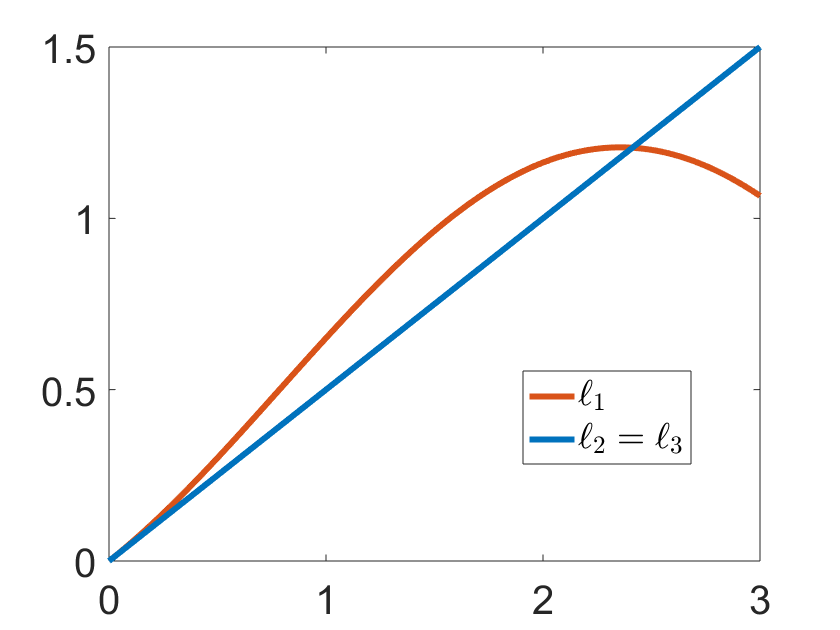}
	\includegraphics[scale=0.28]{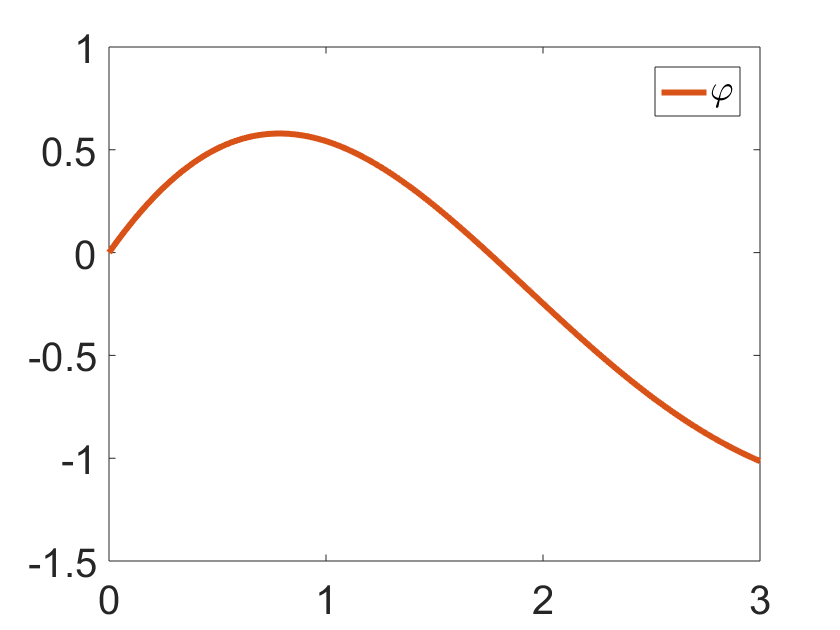}
	\caption{The graphs of control parameters: Example 2
	}
\label{t2}
\end{figure} 
\end{exam}

\begin{exam}
Let us finally discuss the case that all $h_5, h_6, h_8$ are non--zero. Let us choose the constants in  such a way that
$h_1=-{\sqrt{10}\over 4}\cos(t)$, $h_2=-{\sqrt{30}\over 12}\sin(t)+{1 \over 2}$, $h_3=h_4={\sqrt{30}\over 12}\sin(t)+{1 \over 4}$, $h_5={\sqrt{3} \over 3}$, $h_6=h_7=-{\sqrt{3} \over 3}$. We get the solution in the form

\begin{align*}
&x = -{\sqrt{10} \over 4}\sin(t)\\
&y= {\sqrt{30} \over 192}(1-t\sin(t)-\cos(t))+{5\over 64}\cos(t)^2-{5\over 96}\sin(t)\cos(t)-{5t \over 96}+{5\over 48}\sin(t)-{5\over 64} \\
&\theta= -{3\sqrt{10}\over 256}\left((t-8)\sin(t)+\cos(t)-1\right)
\\
& \varphi=
{13\sqrt{3}\over 384}(((t-{96\over 13})\sin(t)+\cos(t)-1)\sqrt{10}\sqrt{3}+({100\over 13}-{50\over 13}\cos(t))\sin(t)-{50\over 13}t)
\\
&\ell_1={\sqrt{30}\over 12}(\cos(t)-1)+{t\over 2} \ \ \ \ 
\ell_2= \ell_3={\sqrt{30}\over 12}(1-\cos(t))+{t \over 4}
\end{align*}

\begin{figure}[ht]
	\includegraphics[scale=0.28]{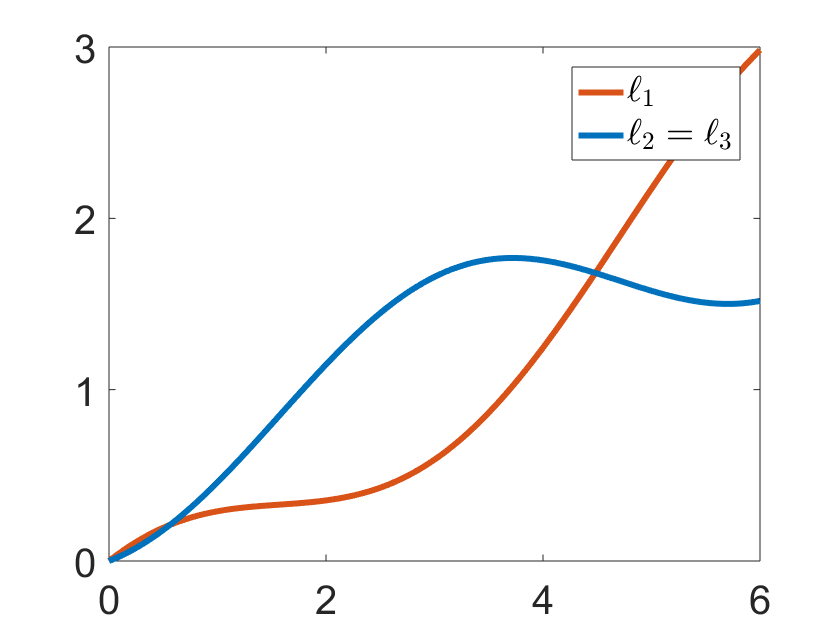}
\includegraphics[scale=0.28]{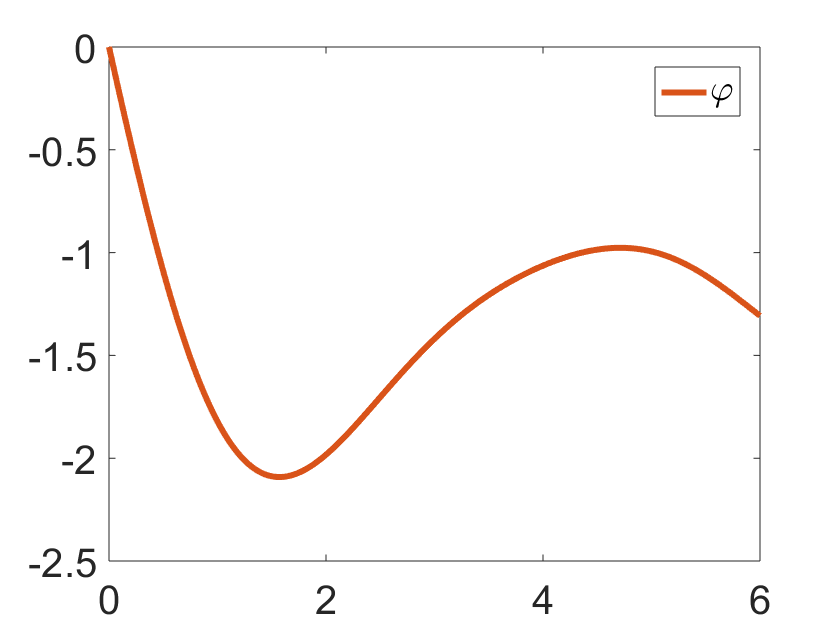}
	\caption{The graphs of control parameters:  Example 3
	}
\label{t3}
\end{figure} 
\end{exam}

\section{Acknowledgment}
The first author is supported by the grant of the Czech Science Foundation no. 17-21360S, `Advances in Snake-like Robot Control' and by a Grant No. FSI-S-17-4464. The second author is partially supported by the grant of the Czech Science Foundation no. 17-01171S, `Invariant differential operators and their applications in geometric modelling and control theory' and by the grant
346300 for IMPAN from the Simons Foundation and the matching 2015-2019 Polish MNiSW fund.

Access to computing and storage facilities owned by parties and projects contributing to the National Grid Infrastructure MetaCentrum provided under the programme "Projects of Large Research, Development, and Innovations Infrastructures" (CESNET LM2015042), is also greatly appreciated.

The algebraic computations are partially calculated in CAS Maple package DifferentialGeometry \cite{dg}.

We would like to thank to Sebastiano Nicolussi Golo and Wojciech Kry\' nski for helpful discussions
and Pawel Nurowski for initial motivation.
We thank the anonymous reviewers whose comments have greatly improved this manuscript.

\end{document}